\documentclass[10pt]{amsart}
\usepackage[active]{srcltx}

\usepackage{amscd, amssymb, amsmath}
\usepackage{epsfig}
\usepackage[all]{xy}

\newtheorem{Df}{Definition}
\newtheorem{theorem}[Df]{Theorem}

\newtheorem{prop}[Df]{Proposition}
\newtheorem{lemma}[Df]{Lemma}

\newtheorem{remark}[Df]{Remark}
\newtheorem{example}[Df]{Example}
\newtheorem{corollary}[Df]{Corollary}

\newcommand{\C}{\ensuremath{\mathbb{C}} }
\newcommand{\Z}{\ensuremath{\mathbb{Z}} }

\newcommand{\N}{\ensuremath{\mathbb{N}} }
\newcommand{\Uhg}{\ensuremath{U_{h}(\g) } }

\newcommand{\slmn}{\ensuremath{\mathfrak{sl}(m|n)}}
\newcommand{\g}{\ensuremath{\mathfrak{g}}}
\newcommand{\slt}{\ensuremath{\mathfrak{sl}(2)}}
\newcommand{\Uslt}{\ensuremath{U_{q}(\slt) } }
\newcommand{\Uqg}{\ensuremath{U_{q}(\g) } }
\newcommand{\UsltH}{\ensuremath{U^H_{q}(\slt) } }

\newcommand{\cat}{\mathcal{C}}
\newcommand{\ideal}{{\overline{A}}}
\newcommand{\links}{\mathcal{L}}
\newcommand{\ob}{Ob(\mathcal{C})}
\newcommand{\obj}{Ob}
\newcommand{\End}{\operatorname{End}}
\newcommand{\Hom}{\operatorname{Hom}}
\newcommand{\unit}{\ensuremath{\mathbb{I}}}
\newcommand{\tr}{\operatorname{tr}}
\newcommand{\Id}{\operatorname{Id}}
\newcommand{\qdim}{\operatorname{qdim}}
\newcommand{\FK}{\mathbb{F}}
\newcommand{\Ubar}{\bar{U}_{q}(\slt)}
\newcommand{\qn}[1]{{\left\{#1\right\}}}
\newcommand{\qd}{\operatorname{\mathsf{d}}}
\newcommand{\osp}{\mathfrak{osp}(2|2n)}
\newcommand{\sll}{\mathfrak{sl}}
\newcommand{\mathsmall}[1]{\mbox{\tiny$#1$}} 
\newcommand{\msmall}[1]{\mbox{\tiny #1}}

\newcommand{\catE}{\mathcal{E}}
\newcommand{\DD}{\mathcal{D}}
\newcommand{\Ho}{{\operatorname{H}}}
\newcommand{\Moc}{{\operatorname{M}}}
\newcommand{\Modd}{{\operatorname{M}'}}
\newcommand{\mc}{{m}}

\newcommand{\lcirc}{\lhd}
\newcommand{\rcirc}{\rhd}
\newcommand{\go}{\rightarrow}

\newcommand{\epsh}[2]
         {\begin{array}{c} \hspace{-1.3mm}
        \raisebox{-4pt}{\epsfig{figure=#1,height=#2}}
        \hspace{-1.9mm}\end{array}}

\newcommand{\pic}[2]{
  \setlength{\unitlength}{#1}
  {\begin{array}{c} \hspace{-1.3mm}
        \raisebox{-4pt}{#2}
        \hspace{-1.9mm}\end{array}}}

\newcommand{\tens}{\otimes}
\newcommand{\rtens}{
\hspace*{0.5ex}
\setlength{\unitlength}{0.08ex}
\begin{picture}(25,18)
 \linethickness{0.2pt}
 \qbezier(9,18)(0,18)(0,9)
 \qbezier(0,9)(0,0)(9,0)
 \qbezier(9,0)(12,0)(24,9)
 \qbezier(9,18)(12,18)(24,9)
 \qbezier(2.636,15.364)(9,9)(15.179,2.821)
 \qbezier(2.636,2.636)(9,9)(15.179,15.179)
\end{picture}
\hspace*{0.5ex}
}
\newcommand{\ltens}{
\hspace*{0.5ex}
\setlength{\unitlength}{0.08ex}
\begin{picture}(25,18)
 \linethickness{0.2pt}
 \qbezier(16,18)(25,18)(25,9)
 \qbezier(25,9)(25,0)(16,0)
 \qbezier(16,0)(13,0)(1,9)
 \qbezier(16,18)(13,18)(1,9)
 \qbezier(22.364,15.364)(16,9)(9.821,2.821)
 \qbezier(22.364,2.636)(16,9)(9.821,15.179)
\end{picture}
\hspace*{0.5ex}
}

\begin{document}
\title{Modified quantum dimensions and re-normalized link invariants}
\author{Nathan Geer}
\address{School of Mathematics\\
Georgia Institute of Technology\\
Atlanta, GA 30332-0160, USA}
\email{geer@math.gatech.edu}
\author{Bertrand Patureau-Mirand}
\address{L.M.A.M., Universit\'e de Bretagne-Sud, BP 573\\
F-56017 Vannes, France }
\email{bertrand.patureau@univ-ubs.fr}
\author{Vladimir Turaev}
\address{IRMA, Universit\'e Louis Pasteur - C.N.R.S.,\\
7 rue Ren\'e Descartes \\
F-67084 Strasbourg, France \\
and \\
Department of Mathematics \\
Indiana University \\
Rawles Hall, 831 East 3rd St \\
Bloomington, IN 47405, USA}

\date{\today}

\begin{abstract}
In this paper we give a re-normalization of the Reshetikhin-Turaev quantum invariants of links, by modified quantum dimensions.  In the case of simple Lie algebras these modified quantum dimensions are proportional to the usual quantum dimensions.  More interestingly we will give two examples where the usual quantum dimensions vanish but the modified quantum dimensions are non-zero and lead to non-trivial link invariants.  The first of these examples is a class of invariants arising from Lie superalgebras previously defined by the first two authors.  These link invariants are multivariable and generalize the multivariable Alexander polynomial.  The second example, is a hierarchy of link invariants arising from nilpotent representations of quantized $\sll(2)$ at a root of unity.  These invariants contain Kashaev's quantum dilogarithm invariants of knots. 
\end{abstract}

\maketitle
\setcounter{tocdepth}{1}



\section{Introduction}
One obstruction to applications of quantum link invariants associated with a
ribbon category $\cat$ stems from the fact that certain simple (irreducible)
objects of $\cat$ may have zero quantum dimensions. If the dimension of a
simple object $V\in \ob$ is zero, then the quantum invariants of all (framed
oriented) links with components labeled by $V$ are equal to zero.  A well
known topological trick allows to derive possibly non-trivial invariants in
this setting, at least in the case of knots. Namely, one presents a
$V$-labeled knot $L$ as the closure of a $(1,1)$ tangle $T$ and considers the
endomorphism of $V$ associated with $T$. This endomorphism is the product of
the identity $\Id_V:V\to V$ with an element $\langle T \rangle$ of the ground
ring of $\cat$.  The tangle $T$ is determined by $L$ uniquely up to isotopy
and therefore $\langle T \rangle$ is an isotopy invariant of $L$. This
invariant may be non-trivial even when $\dim_\cat(V)=0$.  Note that the usual
quantum invariant of $L$ is equal to $\langle T \rangle \dim_\cat(V)$.

For a link $L$ with $\geq 2$ components labeled by $V$, the situation is more
involved because $\langle T \rangle$ may depend on the choice of $T $. In many
known examples of ribbon categories, an appropriate re-normalization of
$\langle T \rangle$ does not depend on the choice of $T $ and yields a
possibly non-trivial invariant of $L$, see \cite{KS}, \cite{ADO}, \cite{Kv},
\cite{GP2}.
A systematic explanation of this phenomenon seems to be missing in the
literature. In this paper we suggest
such an explanation. It is based on a new notion of an ambidextrous object in
$\cat$. Every simple ambidextrous object $J\in \ob$ determines a certain set
$A(J)$ of (isomorphism classes of) simple objects of $\cat$. For all simple
objects $V$ belonging to this set we define a modified (quantum) dimension
depending on $J$. The modified dimension may be non-zero when
$\dim_\cat(V)=0$. Using the modified dimensions we define an isotopy invariant
$F'(L)$ for any link whose components are labeled with objects of $\cat$ under
the only assumption that at least one of the labels belongs to $A(J)$. Most of
these results extend to closed $\cat$-colored ribbon graphs (i.e., to
$\cat$-colored ribbon graphs with no inputs and no outputs).

We give three families of examples illustrating our constructions.  One
example: $\cat$ is the category of finite dimensional $\Uqg$-modules, where
$\Uqg$ be the Drinfeld-Jimbo $\C(q)$-algebra associated to a simple complex
Lie algebra. In this case we recover the standard Reshetikhin-Turaev link
invariants.  In the second example $\cat$ is the category of topologically
free $\Uhg$-modules of finite rank, where $\g$ is a Lie superalgebra of type I
and $\Uhg$ is its quantized universal enveloping $\C[[h]]$-superalgebra.  In
this case we recover the link invariants defined by the first two authors in
\cite{GP2,GP3}. These invariants generalize both the multivariable Alexander
polynomial of links and Kashaev's link invariants.  In the final example
$\cat$ is the category of finite dimensional weight $\Uslt$-modules where $q$
is a root of unity.  We will show that our construction in this case gives a
generalization of the invariants defined by Akutsu, Deguchi and Ohtsuki
\cite{ADO}, using a regularization of the Markov trace and nilpotent
representations of $\Uslt$ at a root of unity.  In the later two examples the
standard Reshetikhin-Turaev link invariant coming from $\cat$ is generically
zero.  

The paper is organized as follows. In Section 2 we recall the basic results
on ribbon categories. In Section 3 we introduce the ambidextrous objects, the
modified dimensions and the invariant $F'$ of closed $\cat$-colored ribbon
graphs. In Section 4 we extend $F'$ to arbitrary $\cat$-colored ribbon graphs
(this does not yield a functor as in the standard theory but only a
quasi-functor). In Section 5 we study the basic properties of $F'$. Section 6
is devoted to the examples.

\section{Ribbon Ab-categories}

We describe the concept of a ribbon Ab-category (for details see \cite{Tu}).
A \emph{tensor category} $\cat$ is a category equipped with a covariant
bifunctor $\otimes :\cat \times \cat\rightarrow \cat$ called the tensor
product, a unit object $\unit$, an associativity constraint, and left and
right unit constraints such that the Triangle and Pentagon Axioms hold.
When the associativity constraint and the left and right unit constraints are
all identities we say the category $\cat$ is a \emph{strict} tensor
category. By Mac Lane's coherence theorem any tensor category is equivalent to
a strict tensor category.

A tensor category $\cat$ is said to be an \emph{Ab-category} if for any pair
of objects $V,W$ of $\cat$ the set of morphisms $\Hom(V,W)$ is an additive
abelian group and the composition and tensor product of morphisms are
bilinear.

Let $\cat$ be a (strict) ribbon Ab-category, i.e. a (strict) tensor
Ab-category with duality, a braiding and a twist. Composition of morphisms
induces a commutative ring structure on $\End(\unit)$.  This ring is called
the \emph{ground ring} of $\cat$ and denoted by $K$.  For any pair of objects
$V,W$ of $\cat$ the abelian group $\Hom(V,W)$ becomes a left $K$-module where
the action is defined by $kf=k\otimes f$ where $k\in K$ and $f\in
\Hom(V,W)$. An object $V$ of $\cat$ is \emph{simple} if $\End(V)= K \Id_V$.

We denote the braiding in $\cat$ by $c_{V,W}:V\otimes W \rightarrow W \otimes
V$ and duality morphisms in $\cat$ by
\begin{align*}
  b_{V} : \:\: & \unit \rightarrow V\otimes V^{*}, & b'_{V} : \:\: &
  \unit\rightarrow V^*\otimes V, & d_{V}: \:\: & V^*\otimes V\rightarrow
  \unit, & d_{V}':\:\: & V\otimes V^{*}\rightarrow \unit.
\end{align*}
The \emph{trace} of any endomorphism $f\in \End(V) $ of an object $V$ of
$\cat$ is defined by
$$\tr_{\cat} (f)= d'_{V} \circ ( f \otimes \Id_{V}^*) \circ b_{V} \in
\End(\unit)= K\, .$$ Define $\dim_{\cat}:\ob \rightarrow K$ by
$\dim_\cat(V)=\tr_{\cat}(\Id_V)$.  We call $\dim_\cat(V)$ the \emph{dimension}
of $V$.

\section{The invariant $F'$ of closed ribbon graphs}\label{SS:inv}
 Let $\cat$ be a strict ribbon Ab-category with ground ring
$K$ and the set of objects $\ob$. We shall assume everywhere that
 $K$ is an integral domain with field of fractions
$\FK$.

For any object $V$ of $\cat$ and any endomorphism $f$ of $V\otimes
V$, set
$$\tr_{L}(f)=(d_{V}\otimes \Id_{V})\circ(\Id_{V^{*}}\otimes
f)\circ(b'_{V}\otimes \Id_{V}) \in \End(V),$$
$$\tr_{R}(f)=(\Id_{V}\otimes d'_{V}) \circ (f \otimes \Id_{V^{*}})
\circ(\Id_{V}\otimes b_{V}) \in \End(V).$$
   An  object $V$ of $\cat$ is called
\emph{ambidextrous} if $\tr_{L}(f)=\tr_R(f)$ for all $f \in\End(V
\otimes V).$

The following lemma gives examples of ambidextrous elements.

\begin{lemma}\label{P:ambdim}
\begin{enumerate}
\item  If $J$ is an   object of $\cat$ such that the braiding
  $c_{J,J}$ commutes with any element of $\End(J \otimes J)$,
  then $J$ is ambidextrous. \label{LI:amb1}
\item If $J$ is a simple object of $\cat$ such that
  $\dim_\cat(J)\neq0$, then $J$ is ambidextrous.  \label{LI:amb2}
\end{enumerate}
\end{lemma}
\begin{proof}

\eqref{LI:amb1} Let $f \in\End(J \otimes J)$.  We have
  $\tr_{R}(f)=\tr_{L}(c_{J,J}^{-1}\circ f \circ c_{J,J}).$
But $c_{J,J}$ commutes with $\End(J\otimes J)$ and so
$c_{J,J}^{-1}\circ f \circ c_{J,J}=f$.

 \eqref{LI:amb2} Let $f \in\End(J \otimes J)$.  We have
  $$\tr_L(f)=\frac{\tr_{\cat}(f)}{\dim_{\cat}(J)}\, \Id_J=\tr_R(f)\, .$$
\end{proof}

Next we recall   the category of $\cat$-colored ribbon graphs
$Rib_\cat$ (for more details see \cite{Tu} Chapter~I).  A morphism
$f:V_1\otimes...\otimes V_n \rightarrow W_1\otimes...\otimes W_m$ in
the category $\cat$ can be represented by the following box and
arrows:
$$\xymatrix{
  \ar@< 8pt>[d]^{W_m}_{... \hspace{1pt}}
  \ar@< -8pt>[d]_{W_1}\\
  *+[F]\txt{ \: \; f \; \;} \ar@< 8pt>[d]^{V_n}_{... \hspace{1pt}}
  \ar@< -8pt>[d]_{V_1}\\
  \: }$$
Such boxes are called coupons. A ribbon graph is formed from several
oriented framed edges colored by objects of $\cat$ and  several
coupons colored with morphisms of $\cat$.   The objects of
$Rib_\cat$  are sequences of pairs $(V,\epsilon)$, where $V\in \ob$
and $\epsilon=\pm$ determines the orientation of the corresponding
edge. The morphism of $Rib_\cat$ are isotopy classes of
$\cat$-colored ribbon graphs and their formal linear combinations
with coefficients in $K$. From now on we write $V$ for $(V,+)$.

Let $F$ be the usual ribbon functor from $Rib_\cat$ to $\cat$ (see
\cite{Tu}). Let $T_{V}$ ($T_V^-$) be a  $\cat$-colored (1,1)-ribbon
graph whose open string is oriented downward (resp. upward) and
colored with a simple object $V$ of $\cat$.  Then $F(T_{V})\in
\End_{\cat}(V)= K \Id_V$.  Let $<T_{V}> \, \in K$ be such that
$F(T_{V})= \, <T_{V}> \, \Id_V$.

Let $V$ and $V'$ be objects of $\cat$ such that $V'$ is simple and define
$$
S'(V,V')=\left< \epsh{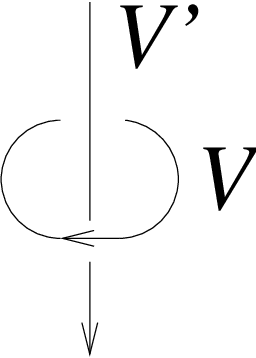}{10ex}\right> \, \in K\, .
$$

\begin{lemma}\label{key}
  For all simple objects $U$, $V$,  $W$  of $\cat$ such that $W$  is
  ambidextrous and for   any $\cat$-colored ribbon graph $T$ with 2
  inputs and 2 outputs colored by $U, V$,
  $$
  S'(U,W)S'(W,V)\left<\, \put(5,20){\msmall{$U$}}
    \hspace{1ex}\epsh{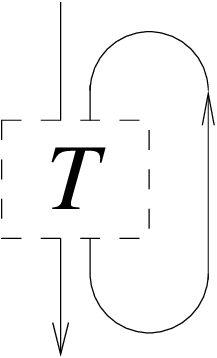}{12ex} \put(1,0){\msmall{$V$}}
    \hspace{2ex}\right> = S'(V,W)S'(W,U) \left< \put(-2,0){\msmall{$U$}}
    \hspace{1ex}\epsh{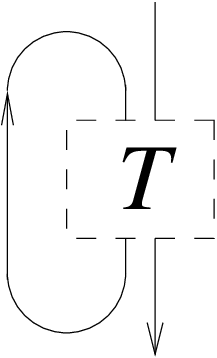}{12ex} \put(-5,20){\msmall{$V$}}
    \hspace{1ex}\right>.
  $$

\end{lemma}
\begin{proof}
  Recall that   $\tr_L(f)=\tr_R(f)$ for all
  $f\in\End(W \otimes W).$ This implies that:
\begin{equation}\label{E:cutV0}
  \left<\, \put(0,15){\msmall{$W$}}
    \hspace{1ex}\epsh{fig37}{8ex} \put(1,0){\msmall{$W$}}
    \put(-10,2){\msmall{$'$}}
    \hspace{2ex}\right> = \left< \put(-3,0){\msmall{$W$}}
    \hspace{2ex}\epsh{fig38}{8ex} \put(-4,15){\msmall{$W$}}
    \put(-3,2){\msmall{$'$}}
    \hspace{1ex}\right>.
\end{equation}
for all $\cat$-colored ribbon graphs $T'$ with 2 inputs and 2 outputs all
colored by $W$.

By definition we have
  \begin{align}
  \label{E:pict1}
  \left<\put(2,15){\msmall{$W$}}
  \put(14,-16){\msmall{$U$}}
  \hspace{2ex}\epsh{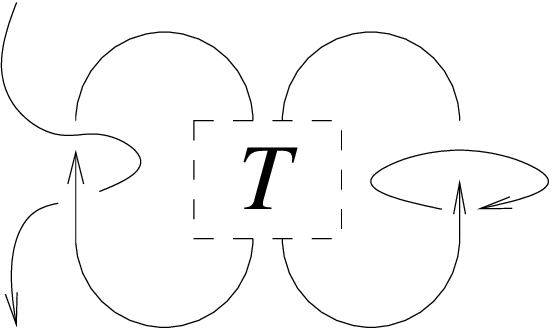}{8ex}
  \put(-10,-14){\msmall{$V$}}
  \put(2,2){\msmall{$W$}}
  \hspace{2ex}\right>&=
  \left< \put(6,15){\msmall{$W$}}
  \hspace{1ex}\epsh{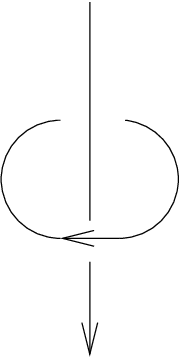}{8ex} \put(2,0){\msmall{$U$}}
  \hspace{2ex}\right>
  \left< \put(4,15){\msmall{$U$}}
  \hspace{1ex}\epsh{fig37}{8ex} \put(2,0){\msmall{$V$}}
  \hspace{2ex}\right>
  \left< \epsh{fig43}{8ex}\put(2,0){\msmall{$W$}}
  \put(-7,15){\msmall{$V$}}
  \hspace{2ex}\right> \notag\\
  &=S'(U,W)\, S'(W,V)\left<\,
  \put(4,15){\msmall{$U$}}
  \hspace{1ex}\epsh{fig37}{8ex} \put(2,0){\msmall{$V$}}
  \hspace{1ex}\right>.
  \end{align}

  Similarly,
  \begin{align}
  \label{E:pict2}
  \left<\put(3,2){\msmall{$W$}}
  \put(15,-14){\msmall{$U$}}
  \hspace{2ex}\epsh{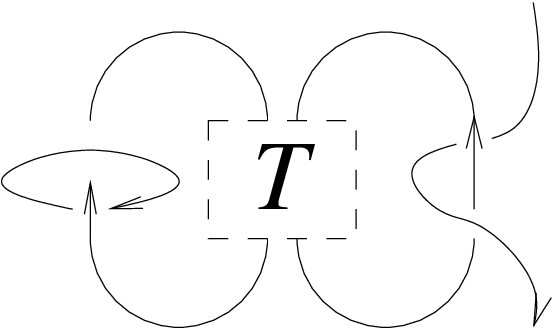}{8ex}
  \put(-10,-14){\msmall{$V$}}
  \put(2,15){\msmall{$W$}}
  \hspace{2ex}\right>  &=
  S'(V,W)\, S'(W,U)\left<
  \put(-2,0){\msmall{$U$}} \hspace{1ex}\epsh{fig38}{8ex}
  \put(-4,15){\msmall{$V$}} \hspace{1ex}\right>.
  \end{align}
  Then equation (\ref{E:cutV0}) implies that the left sides of the above
  equations are equal and so the lemma follows.
\end{proof}

Applying  this lemma to $U=V$, we obtain that if $U$,   $W$ are
simple objects of $\cat$ such that $W$  is
  ambidextrous and  $
  S'(U,W)S'(W,U)\neq 0$, then $U$ is also ambidextrous.

If $A$ is a subset of $\ob$, then let $\links_A$  be the set of
closed $\cat$-colored ribbon graphs, such that at least one of the
colors of the edges is in $A$.  For a simple ambidextrous object $J$
of $\cat$, set  $$ A(J)=\{ V\in \ob : V \text{ simple and } S'(J,V)
\neq 0\}\, .$$ Fix a nonzero   $d_0\in \FK$. For $V\in A(J)$ define
$$\qd_J(V)=d_0\frac{S'(V,J)}{S'(J,V)}\in \FK\, .$$
We view $\qd_J(V)$   as a  modified quantum dimension of $V$
determined by $J$. For any $U,V\in A(J)$, Lemma \ref{key} implies
that
$$ \qd_J(U)\left<\, \put(5,20){\msmall{$U$}}
    \hspace{1ex}\epsh{fig37}{12ex} \put(1,0){\msmall{$V$}}
    \hspace{2ex}\right> = \qd_J(V) \left< \put(-2,0){\msmall{$U$}}
    \hspace{1ex}\epsh{fig38}{12ex} \put(-5,20){\msmall{$V$}}
    \hspace{1ex}\right>
  $$
  for any   $T$.  In
  particular, when $T$ consists of two vertical intervals colored by
  $U, V\in A(J)$,
$$\qd_J(U) \dim_\cat(V) = \qd_J(V) \dim_\cat(U).$$
This shows that the functions $\qd_J$ and  $\dim_\cat$ are
proportional to each other. This is especially interesting when
$\dim_\cat =0$ and $\qd_J\neq 0$, see examples below.

\begin{theorem}\label{T:Main} Let $L\in \links_{A(J)}$ and
  $V$ be the color of an edge of $L$ belonging to ${A(J)}$.  Cutting this
  edge, we obtain
 a colored (1,1)-ribbon graph $T_V$ whose closure is $L$.
  Then $$F'(L) = \qd_J(V)<T_{V}> \in \FK$$ is independent of the choice of the edge 
  to be cut and yields a well defined invariant of $L$.

\end{theorem}
\begin{proof}
  The theorem follows from Lemma \ref{key} and the definition of $\qd_J$.
\end{proof}
We call the invariant $F'$ the re-normalized Reshetikhin-Turaev link invariant.

Let $\ideal=\ideal(J)$ be the set of objects $W$ of $\cat$ such that
there exists a finite family of tuples $(U_i,V_i,f_i,g_i)$, where
$U_i\in  A(J)$, $V_i\in\ob$  and $g_i:U_i \otimes V_i \rightarrow W,
f_i: W\rightarrow U_i \otimes V_i$ satisfying $\Id_W=\sum g_i \circ
f_i$. Note that $\ideal$ has the property that $W\otimes V\in
\ideal$ for all $V\in \ob$ and $W\in \ideal$.

The map $\links_A \to \FK, L\mapsto F'(L)$ extends to a map $\links_\ideal\to
\FK$ as follows.  Let $L$ be a closed $\cat$-colored ribbon graph with one
edge colored by $W\in \ideal$.  Pick a decomposition $\Id_W=\sum g_i \circ
f_i$ as above. Then
$$
\Id_W = F\left(\sum_i \put(2,50){\xymatrix{ \ar[d]^W \\
      *+[F]\txt{ \: \; $g_i$ \; \;} \ar@< 8pt>[d]^{V_i}
      \ar@< -8pt>[d]_{U_i}\\
      *+[F]\txt{ \: \; $f_i$ \; \;}
      \ar[d]^{W}\\
      \: }}\hspace{11ex} \right)\, .
$$
We define $F'(L)$ applying this expansion to the edge of $L$ colored by $W$
and then cutting as above the edge of the resulting graph labeled by $U_i$.
It is easy to show using Theorem \ref{T:Main} that this extension is
independent of the decomposition of $\Id_W$.

\section{The quasi-functor $F'$}
In this section we extend the invariant $F'$ to $\cat$-colored
ribbon graphs with endpoints.  This leads us to a notion of a
quasi-functor which we briefly discuss in a more general setting.
This section is essentially independent from the rest of the paper.

Let $\catE$ be a category.  Given a map
$\Moc:\obj(\catE)\times\obj(\catE)\go\mathcal{S}ets$ we will use
the notation $m\in \Moc$ to mean that there
  exist  two objects $X$ and $Y$ of $\catE$ such that $\mc\in\Moc(X,Y)$.
\begin{Df} \label{D:mod} A $\catE$-bimodule is a map
  $\Moc:\obj(\catE)\times\obj(\catE)\go\mathcal{S}ets$ endowed with two
  operations:
  $$
  \rcirc:\Hom_\catE(Y,Z)\times\Moc(X,Y)\go\Moc(X,Z)
  $$
  $$
  \lcirc:\Moc(Y,Z)\times\Hom_\catE(X,Y)\go\Moc(X,Z)
  $$
  where $X$, $Y$ and $Z$ are any objects of $\catE$.  Given morphisms $f$ and
  $g$ of $\catE$ and $m\in \Moc$ we require that
  \begin{enumerate}
  \item $(f\circ g)\rcirc \mc=f\rcirc(g\rcirc \mc)$.\label{bif1}
  \item $\mc\lcirc (f\circ g)= (\mc\lcirc f)\lcirc g$.\label{bif2}
  \item $(f\rcirc \mc)\lcirc g=f\rcirc( \mc\lcirc g)$.
  \end{enumerate}
  whenever the  operations in these equalities make sense.
\end{Df}
Remark that 
if one has $\Id\rcirc \mc=\mc\lcirc \Id=\mc$ for all $\mc\in\Moc$ these Axioms
mean that $\Moc$ is a bifunctor contravariant in the first place and covariant
in the second place with $\Moc(X,f)(\mc)=f\rcirc \mc$ and
$\Moc(f,Z)(\mc)=\mc\lcirc f$.
\begin{example}
 The functor $\Hom_\catE$ is a $\catE$-bimodule with $\lcirc=\rcirc=\circ$.
\end{example}

\begin{example}\label{Ex:bimod}
  Let $K$ be an integral domain and suppose that $(\Moc, \rcirc, \lcirc)$ is a
  $\catE$-bimodule with values in $K$-modules. If $\Gamma$ is a $K$-module
  then we define an $\catE$-bimodule $\Ho_\catE$ by
  $\Ho_\catE(X,Y)=\Hom_K(\Moc(Y,X),\Gamma)$ with operations $\rcirc$ and
  $\lcirc$ defined as follows.  Let $f:Y\rightarrow Z$ and $f':Z\rightarrow X$
  be morphisms of $\catE$ and let $\phi$ be an element of $\Ho_\catE(X,Y)$
  then
  $$(f\rcirc  \phi)(\mc)=\phi(\mc\lcirc f) \;\; \text{ and } \;\;
  (\phi\lcirc f')(\mc')=\phi(f'\rcirc \mc')$$ where $\mc\in M(Z,X)$ and $\mc'\in
  M(Y,Z)$.
\end{example}

Suppose now that $\catE$ is a tensor $Ab$-category and that $\Moc$
is a $\catE$-bimodule which takes values in   abelian groups. We
assume that the operations $\rcirc$ and $\lcirc$ are bilinear.

\begin{Df}\label{D:monoid}
We call  $\Moc$   a monoidal $\catE$-bimodule if it is endowed with
two bilinear
  operations:
  $$
  \rtens:\Hom_\catE(X,Y)\times\Moc(X',Y')\go\Moc(X\tens X',Y\tens
  Y')
  $$
  $$
  \ltens:\Moc(X,Y)\times\Hom_\catE(X',Y')\go\Moc(X\tens X',Y\tens
  Y')
  $$
  such that for any morphisms $f$, $g$ and $h$ in $\catE$ and any $\mc\in
  \Moc$,
  \begin{enumerate}
  \item $(f\tens g)\rtens\mc=f\rtens (g\rtens\mc)$.
  \item $\mc\ltens(f\tens g)=(\mc\ltens f)\ltens g$.
  \item $f\rtens(\mc\ltens g)=(f\rtens \mc)\ltens g$.
  \item $(f\circ g)\rtens(h\rcirc\mc)=(f\tens h)\rcirc (g\rtens
    \mc)$.\label{ax4}
  \item $(h\rcirc\mc)\ltens(f\circ g)=(h\tens f)\rcirc (\mc\ltens g)$.
  \item $(f\circ g)\rtens(\mc\lcirc h)=(f\rtens \mc)\lcirc (g\tens h)$.
  \item $(\mc\lcirc h)\ltens(f\circ g)=(\mc\ltens f)\lcirc (f\tens g)$.
  \end{enumerate}
whenever the  operations in these equalities make sense.
\end{Df}
  Definition \ref{D:monoid} can be illustrated by pictures.
For example,  Axiom \eqref{ax4} is given by
\begin{center}
  \begin{tabular}{rll}
    \fbox{$\begin{array}{c}f\\\circ\\g\end{array}$}&\rtens&\,
    \fbox{$\begin{array}{c}h\\\nabla\\\mc\end{array}$}\end{tabular}
  $=$\begin{tabular}{c}\fbox{$f\otimes h$}\\$\nabla$\\
    \fbox{$g\rtens \mc$}\end{tabular}
\end{center}
where the composition operations should be read from the top
to the bottom and tensor operations from left to right.
\begin{example}
The functor $\Hom_\catE$ is a monoidal $\catE$-bimodule with
$\lcirc=\rcirc=\circ$ and $\ltens=\rtens=\otimes$.
\end{example}
Suppose that $G:\DD\go\catE$ is a monoidal functor between two tensor
$Ab$-categories.  Let $\Moc$ be a monoidal $\catE$-bimodule and $\Modd$ be a
monoidal $\DD$-bimodule.
\begin{Df}\label{D:quasiF}
  A $G$-bilinear monoidal quasi-functor $G':\Modd\go \Moc$ is a family
  of maps $G':\Modd(X,Y)\go\Moc(G(X),G(Y))$ indexed by the objects $X$
  and $Y$ of $\catE$, such that for every $\mc\in\Modd(X,Y)$ and every
  morphism $f$ of $\DD$ one has:
\begin{enumerate}
  \item $G'(f\rcirc \mc)=G(f)\rcirc G'(\mc)$.
  \item $G'(\mc\lcirc f)=G'(\mc)\lcirc G(f)$.
  \item $G'(f\rtens \mc)=G(f)\rtens G'(\mc)$.
  \item $G'(\mc\ltens f)=G'(\mc)\ltens G(f)$.
  \end{enumerate}
whenever the  operations in these equalities make sense.
\end{Df}

Let us now go  back to the situation of Section \ref{SS:inv}. Define
a  $\cat$-bimodule structure on $\Ho_\cat=\Hom_K(\Hom_\cat,\FK)$ as
in Example \ref{Ex:bimod}.
  In particular, for any objects $U,V$ of $\cat$ we have
$$
\Ho_\cat(U,V)=\Hom_K(\Hom_\cat(V,U),\FK)
$$
and if $f$ and $g$ are morphisms of $\cat$ and $\phi\in\Ho_\cat$ then
$$
(f\rcirc\phi)(g)=\phi(g\circ f)\text{ and } (\phi\lcirc
f)(g)=\phi(f\circ g)
$$
when these operations make sense.

To give $\Ho_\cat$ a monoidal structure, let us recall the partial traces in
$\cat$: if $f\in \Hom_\cat(X\tens Z,Y\tens Z)$ set $\tr_R(f)=(\Id_X\tens
d'_Z)\circ (f\tens \Id_{Z^*})\circ (\Id_X\tens b_{Z})\in\Hom_\cat(X,Y)$ and if
$f\in \Hom_\cat(X\tens Y,X\tens Z)$ set $\tr_L(f)=(d_X\tens\Id_Z )\circ
(\Id_{X^*}\tens f)\circ (b'_{X}\tens\Id_Y)\in\Hom_\cat(Y,Z)$.

Then
$$
\rtens:\Hom_\cat(U,V)\times\Ho_\cat(U',V')\go\Ho_\cat(U\tens U',V\tens V')
$$
$$
\ltens:\Ho_\cat(U,V)\times\Hom_\cat(U',V')\go\Ho_\cat(U\tens U',V\tens V')
$$
are defined as follows.
\begin{enumerate}
\item If $f\in\Hom_\cat(U,V)$ and $\phi\in\Ho_\cat(U',V')$ then $f\rtens\phi
  \in \Ho_\cat(U\tens U',V\tens V')$ is given by
  $$
  (f\rtens\phi)(g)= \phi(\tr_L(g\circ(f\tens \Id_{V'}))).
  $$
  where $g$ is any element of $ \Hom_\cat(V\tens V',U\tens U')$.  This
  operation can be represented by the following diagram:
  $$
  \pic{1.5ex}{
    \begin{picture}(4,4)(2,5)
      \multiput(2,6)(4,0){2}{\line(0,1){2}}
      \multiput(2,6)(0,2){2}{\line(1,0){4}}
      \put(3,6){\vector(0,-1){1}}
      \put(5,6){\vector(0,-1){1}}
      \put(3,8){\line(0,1){1}}
      \put(5,8){\line(0,1){1}}
      \put(3.5,6.8){$g$}
      \put(2,4.5){$\mathsmall{V}$}\put(5.2,4.5){$\mathsmall{V'}$}
      \put(2,9){$\mathsmall{U}$}\put(5.2,9){$\mathsmall{U'}$}
    \end{picture}}
   \mapsto\hspace{1ex}\phi\left(
  \pic{1.5ex}{
    \begin{picture}(5,11)(1,0)
      \put(3,3){\vector(0,-1){1}}
      \put(3,8){\line(0,1){1}}
      \put(3,5){\line(0,1){1}}
      \qbezier(3,2)(3,1)(2,1)
      \qbezier(2,1)(1,1)(1,5.5)
      \qbezier(1,5.5)(1,10)(2,10)
      \qbezier(2,10)(3,10)(3,9)
      \multiput(2,3)(2,0){2}{\line(0,1){2}}
      \multiput(2,3)(0,2){2}{\line(1,0){2}}
      \multiput(2,6)(4,0){2}{\line(0,1){2}}
      \multiput(2,6)(0,2){2}{\line(1,0){4}}
      \put(5,6){\vector(0,-1){6}}
      \put(5,8){\line(0,1){3}}
      \put(2.5,3.7){$f$}\put(3.5,6.8){$g$}
      \end{picture}}
  \right)
  $$
\item If $f\in\Hom_\cat(U',V')$ and $\phi\in\Ho_\cat(U,V)$ then $\phi\ltens f
  \in \Ho_\cat(U\tens U',V\tens V')$ is given by
  $$
  (\phi\ltens f)(g)= \phi(\tr_R(g\circ(\Id_{U'}\tens f)))
  $$
  where $g$ is any element of $\Hom_\cat(V\tens V',U\tens U')$.  Again this
  operation can be represented by the following diagram:
  $$
  \pic{1.5ex}{
    \begin{picture}(4,4)(2,5)
      \multiput(2,6)(4,0){2}{\line(0,1){2}}
      \multiput(2,6)(0,2){2}{\line(1,0){4}}
      \put(3,6){\vector(0,-1){1}}
      \put(5,6){\vector(0,-1){1}}
      \put(3,8){\line(0,1){1}}
      \put(5,8){\line(0,1){1}}
      \put(3.5,6.8){$g$}
      \put(2,4.5){$\mathsmall{V}$}\put(5.2,4.5){$\mathsmall{V'}$}
      \put(2,9){$\mathsmall{U}$}\put(5.2,9){$\mathsmall{U'}$}
    \end{picture}}
   \mapsto\hspace{1ex}\phi\left(
  \pic{1.5ex}{
    \begin{picture}(5,11)(-6,0)
      \put(-3,3){\vector(0,-1){1}}
      \put(-3,8){\line(0,1){1}}
      \put(-3,5){\line(0,1){1}}
      \qbezier(-3,2)(-3,1)(-2,1)
      \qbezier(-2,1)(-1,1)(-1,5.5)
      \qbezier(-1,5.5)(-1,10)(-2,10)
      \qbezier(-2,10)(-3,10)(-3,9)
      \multiput(-2,3)(-2,0){2}{\line(0,1){2}}
      \multiput(-2,3)(0,2){2}{\line(-1,0){2}}
      \multiput(-2,6)(-4,0){2}{\line(0,1){2}}
      \multiput(-2,6)(0,2){2}{\line(-1,0){4}}
      \put(-5,6){\vector(0,-1){6}}
      \put(-5,8){\line(0,1){3}}
      \put(-3.5,3.7){$f$}\put(-4.5,6.8){$g$}
      \end{picture}}
  \right)
  $$ 
\end{enumerate}
One can check that these maps make $\Ho_\cat$ a monoidal $\cat$-bilinear
module.

Fix now a simple ambidextrous object $J$ in $\cat$ and a non-zero
element $d_0$ of $\FK$. Set $\ideal =\ideal (J)$. Let $Rib_\ideal$
be the monoidal $Rib_\cat$-bimodule defined as follows. The
operations of $Rib_\ideal$ are the composition and tensor product of
$Rib_\cat$, i.e. $\lcirc=\rcirc=\circ$ and $\rtens=\ltens=\tens$.
Let $V$ and $W$ be objects of $Rib_\cat$ then $Rib_\ideal(V,W)$ is
the set of $\cat$-colored ribbon graphs in $\Hom_\cat(V,W)$ with at
least one color in $\ideal$.  In particular,
$Rib_\ideal(\emptyset,\emptyset)=\links_\ideal$.
\begin{theorem}
  The invariant $F':\links_\ideal\rightarrow \FK$ extends naturally to a
  $F$-bilinear monoidal quasi-functor $F':Rib_\ideal\go\Ho_\cat$ by the
  formula
  $$
  F'(T)(g)=F'(\tr_{Rib_\cat}(T\circ c_g))
  $$
  where $T\in Rib_\ideal(V,W)$, $g\in \Hom_\cat(F(W),F(V))$ and where $c_g$ is
  a coupon labeled by $g$.  This expresion can be represented by the following
  diagram:
  $$
  F'(T):g\mapsto F'\left(\epsh{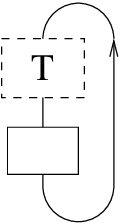}{12ex}\put(-22,-8){g}\right)
  $$
\end{theorem}
\begin{proof}
  Let $T$ and $T'$ be morphisms of $Rib_\cat$ and $Rib_\ideal$, respectively.
  Suppose that $T\circ T'$ is defined.  Let $T''$ be the trace in $Rib_\cat$
  of $T\circ T'$ and let $C$ be a coupon labeled by $F(T')$.  The proof
  follows from the fact that
  $$
  F'(T'')=F'(\tr_{Rib_\cat}(T\circ C))=F'(T)(F(T')).
  $$
  For example, we will now check (1) of Defintion \ref{D:quasiF}.  Let $T$ and
  $T'$ be as above.  Let $g$ be a morphism of $\cat$ and let $c_g$ is a coupon
  labeled by $g$.  Suppose that $g$ is a morphism such that $F'(T\rcirc
  T')(g)$ is defined.  If $c_g$ is a coupon labeled by $g$ then we have
  $$
  \begin{array}{ll}
    F'(T\rcirc T')(g)&= F'(\tr_{Rib_\cat}(T\circ T'\circ c_g))=
    F'(T')(F(c_g\circ T))\\
    &=F'(T')(g\circ F(T))=(F(T)\rcirc F'(T'))(g).
  \end{array}
  $$
  Checking (2)-(4) of Defintion \ref{D:quasiF} are similar.
\end{proof}

\section{Properties of $\qd_J$ and $F'$}

The following lemma shows that the function $\qd_J$ satisfying the conditions of
Theorem \ref{T:Main} is essentially unique.

\begin{lemma}\label{L:d:unique}
  Let $J\in \ob$ be a simple ambidextrous object such that $J\in A(J)$.  Suppose that
  $d:A(J)\rightarrow \FK$ is a function such that the construction of Theorem
  \ref{T:Main} with $\qd_J$ replaced by $d$ yields a well defined invariant
  for all $L\in \links_{A(J)}$. Then $ d=\qd_J $ for an appropriate choice of
  $d_0$.
\end{lemma}
\begin{proof} Let $L$ be the Hopf link with components colored $V $ and $J$,
  where $V$ is a simple object of $\cat$.  By opening one strand of $L$ and
  then the other, we get
  $$
  d(V) \, S'(J,V)=d(J)\, S'(V,J)\, .
  $$
  So, $d=d_0^{-1}d(J)\qd_J$.
\end{proof}

This lemma implies in particular that if $J_1, J_2$ are two simple
ambidextrous objects such that $A(J_1)=A(J_2)$, then $\qd_{J_1}=\qd_{J_2}$ for
any choice of $d_0$ in the definition of $\qd_{J_1}$ and an appropriate choice
of $d_0$ in the definition of $\qd_{J_2}$.

\begin{lemma}\label{L:duality}
  Let $V,W$ be simple objects in $\cat$, then $S'(V,W)=S'(V^*,W^*)$.
\end{lemma}
\begin{proof}
  We have
  $$
  S'(V,W)=\left< \put(4,15){\msmall{$W$}} \hspace{1ex}\epsh{fig43}{8ex}
    \put(1,0){\msmall{$V$}} \hspace{2ex}\right> = \left<
    \put(4,15){\msmall{$W$}} \hspace{1ex}\epsh{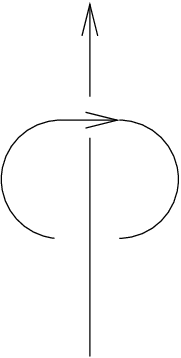}{8ex}
    \put(1,0){\msmall{$V$}} \hspace{2ex}\right> = \left<
    \put(0,16){\msmall{$W^*$}} \hspace{1ex}\epsh{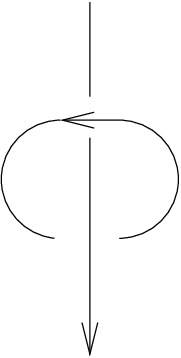}{8ex}
    \put(1,0){\msmall{$V^*$}} \hspace{2ex}\right> = S'(V^*,W^*)
  $$
  where the first and the fourth equalities follow from the definition of
  $S'$, the second from composing with the morphism $(\Id_V \otimes d_V)(b_V
  \otimes \Id_V)$ and the third from the property that
  $F(\uparrow_V)=F(\downarrow_{V^*})$.
\end{proof}
\begin{lemma}\label{L:dualityBIS}
  If an object $J$ in $\cat$ is ambidextrous, then so is $J^*$.
\end{lemma}
\begin{proof} We have to prove that $\tr_L(f)=\tr_R(f)$ for any endomorphism
  $f$ of $J^*\otimes J^*$.  First we define an endomorphism $f'$ of
  $J\otimes J$ by the following formula:
  $$f'=F\left(
  \pic{2ex}{
    \begin{picture}(8,6)(-1,0)
      \multiput(2,2)(2,0){2}{\line(0,1){2}}
      \multiput(2,2)(0,2){2}{\line(1,0){2}}
      \qbezier(2.5,2)(2.5,1.25)(2,1.25)
      \qbezier(2,1.25)(1,1.25)(1,5)
      \put(1,6){\vector(0,-1){1}}
      \qbezier(3.5,2)(3.5,0.5)(2,0.5)
      \qbezier(2,0.5)(0,0.5)(0,5)
      \put(0,6){\vector(0,-1){1}}
      \qbezier(2.5,4)(2.5,5.5)(4,5.5)
      \qbezier(4,5.5)(6,5.5)(6,1)
      \put(6,1){\vector(0,-1){1}}
      \qbezier(3.5,4)(3.5,4.75)(4,4.75)
      \qbezier(4,4.75)(5,4.75)(5,1)
      \put(5,1){\vector(0,-1){1}}
      \put(2.7,2.7){$f$}
      \put(-1,5){$J$}\put(1.2,5){$J$}
      \put(4,0){$J$}\put(6.2,0){$J$}
    \end{picture}
  }\right).
  $$
  As $V$ is ambidextrous we have $\tr_L(f')=\tr_R(f')$ and thus
  $$F\left(\pic{2ex}{
    \begin{picture}(6,6)(1,0)
      \multiput(2,2)(2,0){2}{\line(0,1){2}}
      \multiput(2,2)(0,2){2}{\line(1,0){2}}
      \qbezier(2.5,2)(2.5,1.25)(2,1.25)
      \qbezier(2,1.25)(1,1.25)(1,5)
      \put(1,6){\vector(0,-1){1}}
      \qbezier(2.5,4)(2.5,5.5)(4,5.5)
      \qbezier(4,5.5)(6,5.5)(6,1)
      \put(6,1){\vector(0,-1){1}}
      \qbezier(3.5,4)(3.5,4.75)(4,4.75)
      \qbezier(4,4.75)(5,4.75)(5,3)
       \qbezier(5,3)(5,1.25)(4,1.25)
      \qbezier(4,1.25)(3.5,1.25)(3.5,2)
      \put(5,3.05){\vector(0,-1){0.1}}
      \put(2.7,2.7){$f$}\put(1.2,5){$J$}\put(6.2,0){$J$}
    \end{picture}
  }\right)=F\left(\pic{2ex}{
    \begin{picture}(6,6)(-1,0)
      \multiput(2,2)(2,0){2}{\line(0,1){2}}
      \multiput(2,2)(0,2){2}{\line(1,0){2}}
      \qbezier(2.5,2)(2.5,1.25)(2,1.25)
      \qbezier(2,1.25)(1,1.25)(1,3)
      \qbezier(1,3)(1,4.75)(2,4.75)
      \qbezier(2,4.75)(2.5,4.75)(2.5,4)
      \qbezier(3.5,2)(3.5,0.5)(2,0.5)
      \qbezier(2,0.5)(0,0.5)(0,5)
      \put(0,6){\vector(0,-1){1}}
      \qbezier(3.5,4)(3.5,4.75)(4,4.75)
      \qbezier(4,4.75)(5,4.75)(5,1)
      \put(5,1){\vector(0,-1){1}}
      \put(1,3.05){\vector(0,-1){0.1}}
      \put(2.7,2.7){$f$}\put(-1,5){$J$}\put(4,0){$J$}
    \end{picture}
  }\right)$$
In ribbon categories, modules are canonically isomorphic to their bidual and 
if $T$ is a colored ribbon graph, changing both the orientation of an internal edge
and its color to its dual does not affect $F(T)$.
Doing this for the internal edges of the two ribbon graphs above, we deduce that
$$
\left<\pic{2ex}{
    \begin{picture}(5,5)(1,0.5)
      \multiput(2,2)(2,0){2}{\line(0,1){2}}
      \multiput(2,2)(0,2){2}{\line(1,0){2}}
      \put(2.5,1){\vector(0,1){1}}
      \put(2.5,4){\vector(0,1){1}}
      \qbezier(3.5,4)(3.5,4.75)(4,4.75)
      \qbezier(4,4.75)(5,4.75)(5,3)
      \qbezier(5,3)(5,1.25)(4,1.25)
      \qbezier(4,1.25)(3.5,1.25)(3.5,2)
      \put(5,2.95){\vector(0,1){0.1}}
      \put(2.7,2.7){$f$}\put(5.3,2.7){$J^*$}
      \put(1.2,5){$J$}\put(1.2,0.5){$J$}
    \end{picture}
  }\right>=\left<\pic{2ex}{
    \begin{picture}(5,5)(0,0.5)
      \multiput(2,2)(2,0){2}{\line(0,1){2}}
      \multiput(2,2)(0,2){2}{\line(1,0){2}}
      \qbezier(2.5,2)(2.5,1.25)(2,1.25)
      \qbezier(2,1.25)(1,1.25)(1,3)
      \qbezier(1,3)(1,4.75)(2,4.75)
      \qbezier(2,4.75)(2.5,4.75)(2.5,4)
      \put(1,3.05){\vector(0,-1){0.1}}
      \put(3.5,1){\vector(0,1){1}}
      \put(3.5,4){\vector(0,1){1}}
      \put(2.7,2.7){$f$}\put(-0.5,2.7){$J^*$}
      \put(3.8,5){$J$}\put(3.8,0.5){$J$}
    \end{picture}
  }\right>
$$
Which mean that $\tr_R(f)=\tr_L(f)$.
\end{proof}
\begin{corollary}\label{bbb} For any   simple ambidextrous object $J\in \ob$,
  we have $A (J^*)=\{ V^* : V \in A(J) \}$ and $\qd_{J^*}(V^*)=\qd_J (V)$ for
  all $V\in A (J)$. If $V, V^*\in A(J)$ satisfy $S'(V,V^*)\neq 0$, then
  $\qd_J(V^*)=\qd_J (V)$.
\end{corollary}
\begin{proof} Lemma~\ref{L:duality} implies that $$A (J^*)= \{ V\in \ob : V
  \text{ simple and } S'(J^*,V) \neq 0\}= \{ V^* : V \in A(J) \}$$ and for all
  $V\in A (J)$,
  $$
  \qd_J (V)=d_0\frac{S'(V,J)}{S'(J,V)}=
  d_0\frac{S'(V^*,J^*)}{S'(J^*,V^*)}= \qd_{J^*}(V^*)\, .
  $$

  Consider the Hopf link $H$ with components labeled by $V$ and $V^*$.  Then
  by definition, $F'(H)$ is equal to both $\qd_J (V)S'(V^*,V)$ and $\qd_J
  (V^*)S'(V,V^*)$.  Now Lemma~\ref{L:duality} implies that
  $S'(V^*,V)=S'(V,V^*)\neq 0$ and therefore $\qd_J(V^*)=\qd_J (V)$.
\end{proof}
Let $I$ be the set of isomorphism classes of simple objects of $\cat$. We call
a subset $B$ of $I$ \emph{complete} if $S'(U,V)\neq 0$ for all $U,V\in B$ and
$S'(U,W)= 0$ for all $U\in B$, $W\in I-B$.
\begin{lemma}\label{L:B}
  Let $B$ be a complete subset of the set $I$, which contains at least one
  ambidextrous object. Then
  \begin{enumerate}
  \item all objects in $B$ are ambidextrous; \label{LI:B1}
  \item for any object $J\in B$, we have $A(J)=B$; \label{LI:B2}
  \item Let $F'_J, F'_{J'}$ be the invariants derived above from arbitrary
    pairs $(J \in B, d_0)$ and $(J' \in B, d'_0)$, respectively, where $d_0,
    d'_0\in \FK\setminus\{0\}$.  Then
    $$F'_J=\frac{\qd_J(J')}{d_0'}\, F'_{J'}\, .$$\label{LI:B3}
  \end{enumerate}
\end{lemma}
\begin{proof}
  The first claim follows from the definition of a complete set and the remark
  following Lemma \ref{key}. The equality $A(J)=B$ follows from the
  definitions.

  Lemma \ref{L:d:unique} implies that $\qd_{J'}$ is proportional to $\qd_J$.
  More precisely, $\qd_J=(d_0')^{-1}\qd_J (J')\qd_{J'}$.
  Thus,$$F'_J(L)=\qd_{J}(V)<T_V>= (d_0')^{-1}\qd_J
  (J')\qd_{J'}(V)<T_V>=(d_0')^{-1}\qd_J (J')F'_{J'}(L).$$
\end{proof}

The completeness condition on $B$ seems to be very strong. However, complete sets
naturally arise in our examples, see Subsections \ref{SS:E1}, \ref{SS:E2} and
\ref{SS:E3}.  In such a situation Lemma~\ref{L:B} says that the ambidextrous
object $J$ defines a ``cluster'' $A(J)$ of ambidextrous objects with the
property that any element of this cluster leads to an invariant proportional
to $F'_J$. There are two kinds of clusters depending on whether or not
$\dim_\cat(J)=0$. The case $\dim_\cat(J)=0$ is of particular interest because
then $F'$ may be non-zero while the usual invariant $F$ restricted to
$\links_\ideal$ is zero, as will be clear from the next lemma.

From now on and up to the end of this section we fix a simple ambidextrous
object $J$ in $\cat$.

\begin{lemma}\label{P:dim0}  For any $V\in A(J)$,
$$\dim_\cat(V)=d_0^{-1}\dim_\cat(J)\,\qd_J (V)\, .$$    
\end{lemma}
\begin{proof}
  Consider the Hopf link $H$ with
  components labeled by $J$ and $V\in A$.  Now $F(H)$ can be computed in two
  ways, namely by cutting the component labeled by $J$ or cutting the
  component labeled by $V$.  This gives
  $$F(H)=S'(V,J)\dim_\cat(J)= S'(J,V)\dim_\cat(V)\, .$$  Thus,
  $\dim_\cat(V)=d_0^{-1}\dim_\cat(J)\,\qd_J(V)$.
  \end{proof}
The corollary below follows directly from the previous lemma.
\begin{corollary} \label{C:dim0}
The following hold,
  \begin{enumerate}
  \item if $\dim_\cat(J)=0$ then $\dim_\cat(V)=0$ for all $V\in A(J)$ and
    $F(L)=0$ for all $L$ in $\links_\ideal$,
  \item if $\dim_\cat(J)\neq 0$ then $F'$ is proportional to $F$.
  \end{enumerate}
 \end{corollary}

The following two propositions show that $F'$ has behavior similar to the
functor $F$.

\begin{prop}\label{P:cab}
  Let $U,V,W\in \ob$ be such that $W\cong U\otimes V$.  Let $L$ be a
  $\cat$-colored link, such that a component of $L$ is colored by $W$.  Let
  $L_\parallel$ be the link obtained from $L$ by replacing this component of $L$ by
  two parallel copies colored by $U$ and $V$.  If $L_\parallel$ is an element of
  $\links_{\ideal(J)}$, then $F'(L)=F'(L_\parallel)$.
\end{prop}
\begin{proof}
  Since $W\cong U\otimes V$ there exist morphisms $f:W\rightarrow U\otimes V$
  and $g: U\otimes V\rightarrow W$ such that $f\circ g=\Id_{U\otimes V}$ and
  $g\circ f=\Id_W$.  Use the equality $g\circ f=\Id_W$ to replace a portion of
  the $W$-colored component of $L$ by two strings labeled by $U$ and $V$ and
  two coupons labeled $f$ and $g$.  Then by sliding one of the coupons around
  the component and using the equality $f\circ g=\Id_{U\otimes V}$ one arrives
  at $L_\parallel$.  Then $F'(L)=F'(L_\parallel)$ since $F'$ is a well defined $\cat$-colored
  ribbon graph invariant.
\end{proof}
\begin{prop}\label{P:S'non}
  Let $L$ be an element of $\links_{A(J)}$ with a circle component colored by
  $W\in \ob$.  If $W$ is the only color of $L$ in $A(J)$, then we additionally
  assume that $W^*\in A(J)$ and $S'(W, W^*)\neq 0$. Let $L_-$ be obtained from
  $L$ by reversing the orientation of the $W$-colored component and changing
  its color to $W^*$. Then $F'(L)=F'(L_-)$.
\end{prop}
\begin{proof} We consider two cases: (a) if $W$ is the only color of $L$ in
  $A(J)$ and (b) otherwise.  In the latter case we have
  $F'(L)=\qd_J(V)\left<F(T_V)\right>$ where $V\in A(J)$ is the label of
  another circle component (or an edge) of $L$ and $T_V$ is obtained from $L$
  by cutting this component (edge).  But $F(T_V)=F(T'_V)$ where $T'_V$ is
  obtained from $T_V$ by reversing the orientation of the $W$-colored
  component and changing its color to $W^*$.  In the former case we have
  \begin{align*}
    F'(L)&=\qd_J(W)\left<T_W\right>\\
    &=\qd_J(W)\left<T_{W^*}^-\right>\\
    &=\qd_J(W^*)\left<\bot_{W^*}\right>\\
    &= F'(L_-)
  \end{align*}
  where $\bot_{W^*}$ is the ribbon graph $T_{W^*}^-$ rotated $180^\circ$. (We
  use here the second claim of Corollary \ref{bbb}).
\end{proof}

\begin{prop} If $ L\in \links_{\ideal(J)}$, $L_+\in \links_{\ob}$, then the
  disjoint union $L\sqcup L_+$ belongs to $\links_{\ideal (J)}$ and $F'(L
  \sqcup L_+)=F'(L)\, F(L_+)$.
\end{prop}
\begin{proof}
  Follows from the definitions.
\end{proof}

\begin{remark}
  Both Propositions \ref{P:cab} and \ref{P:S'non} can be extended to analogous
  statements for non-closed $\cat$-colored ribbon graphs. In other words, $F'$
  behaves under cabling and reversing orientation in the same way as the
  standard ribbon functor $F$.
\end{remark}

\section{Examples}
We give three classes of examples of ambidextrous objects and associated
re-normalized link invariants.

\subsection{Link invariants   from Lie algebras}\label{SS:E1}
Let $\g$ be a simple Lie algebra and let $\Uqg$ be the Drinfeld-Jimbo
$\C(q)$-algebra associated to $\g$ (see \cite{Tu}, XI.6).  
(Note here $q$ is not a root of unity.)  Let $\cat$ be the
category of finite dimensional $\Uqg$-modules.  It is well known that $\cat$
is a ribbon Ab-category with ground ring $K=\C(q)$.  Here $\FK=K$ and
$\tr_\cat$ (resp. $ \dim_\cat)$ is the quantum trace (resp.\ quantum
dimension).
\begin{lemma}
  All simple objects of $\cat$ are ambidextrous.
\end{lemma}
\begin{proof}
  This follows from Lemma \ref{P:ambdim}, as $\dim_\cat(J)\neq 0$ for any
  simple object $J$ of $\cat$.
\end{proof}
Let $I$ be the set of isomorphism classes of simple objects of $\cat$.  One
can show that $S'(V,W)\neq 0$ for any $V, W \in I$ (for a similar calculation
see \cite{GP2}, Proposition 2.2). Thus, $ A(J)=I$ for all $J\in I$. The
construction of Section~\ref{SS:inv} derives from $J$ and any non-zero $d_0\in
\FK$ a function $\qd_J$ and an invariant $F'$. Lemma~\ref{L:B} implies that
$F'$ is essentially independent of the choice of $J$.
\begin{prop}
  For any $J\in I$ and $d_0=\qdim(J)$, the $J$-determined quantum dimension
  $\qd_J$ is equal to the usual quantum dimension and $F'$ is the usual
  Reshetikhin-Turaev quantum group invariant arising from $\g$.
\end{prop}
\begin{proof}
  This follow from Corollary \ref{C:dim0}.
\end{proof}

\subsection{Link invariants   from Lie superalgebras}\label{SS:E2}
In \cite{GP2,GP3} the first two authors derived new link invariants from Lie
superalgebras $\osp$ and $\slmn$ with $m\neq n$. In particular, the invariants
associated with $\sll(m|1)$ generalize both the multivariable Alexander
polynomial of links and Kashaev's link invariants. We explain here that the
construction of \cite{GP2,GP3} is a special case of the construction of
Section \ref{SS:inv}.

In this subsection we work in the category of vector superspaces with even
morphisms, i.e., the category whose objects are $\Z/2\Z$-graded vector spaces
and the morphisms are even linear maps. A Lie superalgebra $\g=\slmn$ or
$\g=\osp$ gives rise to the quantized universal enveloping
$\C[[h]]$-superalgebra $\Uhg$. Let $\cat$ be the category of topologically
free $\Uhg$-modules of finite rank (i.e., modules of the form $V[[h]]$ where
$V$ is a finite dimensional $\g$-module). It is known that $\cat$ is a ribbon
Ab-category with ground ring $K=\C[[h]]$, see for instance \cite{GP2} and
references there. The object $V[[h]]$ of $\cat$ is simple if and only if $V$
is a simple $\g$-module.  The finite dimensional simple $\g$-modules are
divided into two classes: typical and atypical. A simple $\g$-module is
typical if each time it is a submodule or a factor-module of a finite
dimensional $\g$-module, it splits as a direct summand. We call a
$\Uhg$-module $V[[h]]$ (a)typical if $V$ is a (a)typical $\g$-module.

Let $I$ be the set of isomorphism classes of simple objects of $\cat$ and let
$B$ be the subset of $ I$ consisting of isomorphism classes of typical
$\Uhg$-modules.
\begin{lemma}
  If $V$ is an element of $B$, then $\dim_\cat(V)=0$. The link invariant $F$
  restricted to $\links_B$ is zero.
\end{lemma}
\begin{proof}
  The first statement follows from a direct calculation using the character
  formula of $V$ (see \cite{GP2}).  The second statement follows from Corollary 
  \ref{C:dim0}.
\end{proof}
The following lemma is a restatement of Lemma 2.8 in \cite{GP2}.
\begin{lemma}\label{L:VnoMulti}
  There exists an element $J_0\in B$ such that $J_0\otimes J_0$ splits as a
  direct sum of elements of $B$ with multiplicity one.  In particular, the
  algebra $\End_{\cat}(J_0\otimes J_0)$ is commutative.
\end{lemma}
\begin{corollary}
  The object $J_0$ of $\cat$ is ambidextrous.
\end{corollary}
\begin{proof}
  Lemma \ref{L:VnoMulti} implies that the braiding $c_{J_0,J_0}$ commutes with
  all elements of $\End_{\cat}(J_0\otimes J_0)$.  Thus, the corollary follows
  from Lemma~\ref{P:ambdim}.
\end{proof}
In \cite{GP2} the first two authors have shown that $S'(U,V)\neq 0$ for all
$U,V\in B$ and $S'(U,W)= 0$ for all $U\in B$, $W\in I-B$. In other words, the
set $B$ is complete.
\begin{prop}
  Every element $J\in B$ is ambidextrous and $ A(J)=B$.  The construction in
  Section~\ref{SS:inv} gives a function $\qd=\qd_J:B\rightarrow
  \C[[h]][h^{-1}]$ and an invariant $F'=F'_J$, which do not depend on $J$ up
  to multiplication by a non-zero element of $\C[[h]][h^{-1}]$.
\end{prop}
\begin{proof}
  Since $B$ is complete and contains $J$ the lemma follows from
  Lemma~\ref{L:B}.
\end{proof}

The link invariant introduced in \cite{GP2} is defined on $\links_B$. Its
definition is similar to the one of $F'$ above but uses a certain function
$d:B\rightarrow \C[[h]][h^{-1}]$ rather than $\qd$.  By Lemma
\ref{L:d:unique}, the function $d$ must be proportional to $\qd$ and therefore
the link invariant of \cite{GP2} is equal to the invariant $F'$ associated
with an arbitrary $J\in B$ and an appropriate $d_0$ (depending on $J$).

\begin{remark}
  (1) In this example there is no canonical choice for $J\in B$.  However, for
  each $J$ there is a suitable choice of $d_0$ (possibly, distinct from
  $d_0=1$) such that $\qd_J$ has a nice formula (cf.\ \cite{GP2}). This
  justifies our choice to include the factor $d_0$
  in the definition of $F'$.\\
  (2) The extension of $F'$ in Section \ref{SS:inv} could be useful in
  computing   $F'$ for links colored with non semisimple modules. \\
  (3) It would be interesting to extend the constructions of this subsection
  to other Lie superalgebras.
\end{remark}

\subsection{Link invariants from $\Uslt$ at roots of unity}\label{SS:E3}
In this subsection we will consider the generalized multivariable Alexander
invariants defined by Akutsu, Deguchi and Ohtsuki in \cite{ADO}, which
contain Kashaev's invariants (see \cite{Kv,MM}).  These invariants are
indexed by positive integers.  In \cite{jM}, Jun Murakami gives a framed
version of these invariants using the universal $R$-matrix of $\Uslt$ and
calls them the colored Alexander invariants.  Here we show that these
invariants are restrictions of invariants defined using ribbon categories as
formulated above.

Fix a positive integer $N$ and let $q=e^\frac{\pi\sqrt{-1}}{N}$ be a $2N$-root
of unity.  We use the notation $q^x=e^{\frac{\pi\sqrt{-1} x}{N}}$.  Here we
give a slightly generalized version of $\Uslt$.  Let $\UsltH$ be the
$\C(q)$-algebra given by generators $E, F, K, K^{-1}, H$ and relations:
\begin{align*}
  HK&=KH, & HK^{-1}&=K^{-1}H, & [H,E]&=2E, & [H,F]&=-2F,\\
  KK^{-1}&=K^{-1}K=1, & KEK^{-1}&=q^2E, & KFK^{-1}&=q^{-2}F, &
  [E,F]&=\frac{K-K^{-1}}{q-q^{-1}}.
\end{align*}
The algebra $\UsltH$ is a Hopf algebra where the coproduct, counit and
antipode are defined by
\begin{align*}
  \Delta(E)&= 1\otimes E + E\otimes K, & \Delta(F)&=K^{-1} \otimes F +
  F\otimes 1,\\
  \Delta(H)&=H\otimes 1 + 1 \otimes H, & \Delta(K)&=K\otimes K, &
  \Delta(K^{-1})&=K^{-1}\otimes K^{-1},\\
  \epsilon(E)&= \epsilon(F)=\epsilon(H)=0, &
  \epsilon(K)&=\epsilon(K^{-1})=1,\\
  S(E)&=-EK^{-1}, & S(F)&=-KF, & S(K)&=K^{-1}.
\end{align*}
Define $\Ubar$ to be the Hopf algebra $\UsltH$ modulo the relations
$E^N=F^N=0$.

We say a $\Ubar$-module $V$ is a \emph{weight module} if $V$ has a weight
decomposition with respect to $H$ and if $q^H$ acts as $K$.  Let $\cat$ be the
tensor Ab-category of finite dimensional weight $\Ubar$-modules (here the
ground ring is $\C$).  We say a simple weight module is \emph{typical} if its
highest weight is in the set $(\C \setminus \Z) \cup \{-1 + kN : k\in \Z\}$
otherwise we say it is \emph{atypical}.  A typical module is $N$ dimensional
and indexed by its highest weight $\lambda$.  We denote such a module by
$V_\lambda$ (for a basis of this module see \cite{jM}).  The weights of
$V_\lambda$ are $\lambda -2i$ for $0 \leq i \leq N-1$, so its character
formula is $\sum_{i=0}^{N-1}u^{\lambda -2i}$ where the coefficient of $u^a$ is
the dimension of the $a$-weight space.

We will now recall that the category $\cat$ is a ribbon Ab-category.  
For $a\in \C$ and $n\in\N$ we set $\{a\}=q^a-q^{-a}$ and $\{n\}!=\{n\}\{n-1\}...\{1\}$.
In \cite{O} Ohtsuki defines an element $R_t$ given by
$$
R_t=q^{H\otimes H/2} \sum_{n=0}^{N-1} \frac{\{1\}^{2n}}{\{n\}!}q^{n(n-1)/2}
E^n\otimes F^n,
$$
where $q^{H\otimes H/2}$ is a formal symbol.  If $v$ and $v'$ are two weight vectors of weights  of $\lambda$ and $\lambda'$ then $q^{H\otimes H/2}$ acts on $v\otimes v'$ by 
$$q^{H\otimes H/2}.(v\otimes v') =q^{\lambda \lambda'/2}v\otimes v'.$$
 Thus, the action of $R_t$ on the tensor product of two objects of $\cat$ is well
defined and induces an endomorphism on such a tensor product.  Moreover, $R_t$    
gives rise to a braiding $c_{V,W}:V\otimes W \rightarrow
W \otimes V$ on $\cat$ defined by $v\otimes w \mapsto \tau(R_t(v\otimes w))$
where $\tau$ is the permutation $x\otimes y\mapsto y\otimes x$ (see \cite{jM,O}).   

\begin{remark}
  It is important that $\Ubar$ contains $H$ because the modules $V_\lambda$
  and $V_{\lambda +2N}$ over $\Uslt/\{E^N=F^N=0\}$ are isomorphic but the
  action of the $R$-matrix on $ V_\lambda^{\otimes 2}$ and $V_{\lambda
    +2N}^{\otimes 2}$ are different.  These modules are distinct in $\cat$ as
  $H$ acts differently.
\end{remark}

Let $V$ and $W$ be objects of $\cat$.  Let $\{v_i\}$ be a basis of $V$ and
$\{v_i^*\}$ be a dual basis of $V^*$.  Then
\begin{align*}
  b_{V} :& \C \rightarrow V\otimes V^{*}, \text{ given by } 1 \mapsto \sum
  v_i\otimes v_i^* & d_{V}: & V^*\otimes V\rightarrow \C, \text{ given by }
  f\otimes w \mapsto f(w)
\end{align*}
are duality morphisms of $\cat$.  

Also, in \cite{O} Ohtsuki defines an element $u$ given by
$$
\sum_{n=0}^{N-1} S(t_n)q^{-H^2/2}s_n
$$
where $R_t=q^{H\otimes H/2} \sum_{n=0}^{N-1} s_n \otimes t_n$ and $q^{-H/2}$ is a formal symbol whose action on a weight vector $v_\lambda$ is given by $q^{-H^2/2}.v_\lambda = q^{-\lambda^2/2}v_\lambda.$
Let $\theta=uK^{N-1}=K^{N-1}u$. 
The twist $\theta_V:V\rightarrow V$ is defined by $v\mapsto \theta^{-1}v$ (see \cite{jM,O}). 

If $L$ is a link colored with objects of $\cat$ such that one of the colors is
typical then $F(L)=0$ where as above $F$ is the usual ribbon functor
$F:Rib_\cat\rightarrow\cat$.  We will now show that the general construction
above gives rise to a non-trivial invariant $F'$, which contains the ADO
invariants and so Kashaev's invariants.

Fix a typical $\Ubar$-module $V_{\lambda_0}$ such that $\lambda_0 \in
\C\setminus \frac12\Z$ and denote it by $J_0$.
\begin{lemma}\label{L:VnoMult3b}
  The tensor product $J_0\otimes J_0$ splits as a direct sum of typical
  $\Ubar$-modules with no multiplicity.  In particular, the algebra
  $\End_{\cat}(J_0\otimes J_0)$ is commutative.
\end{lemma}
\begin{proof}
  First, using the character formula for $J_0$ one sees that all the weights
  of $J_0\otimes J_0$ are not integral.  Then since typical modules always
  split we have that $J_0\otimes J_0$ is a direct sum of typical modules.  The
  character formula for a typical module then implies that $J_0\otimes
  J_0=\oplus_{i=0}^{N-1} V_{2\lambda_0 - 2i}$, which completes the proof.
\end{proof}
\begin{corollary}
  The element $J_0$ is ambidextrous.
\end{corollary}
\begin{proof}
  The corollary follows directly from Lemma \ref{P:ambdim} and
  Lemma~\ref{L:VnoMult3b}.
\end{proof}
Next we will compute $dim_\cat(V_\lambda)$ and
$S'(V_\lambda,V_{\lambda'})$.  
To do this we will need the morphisms $d_{V}': V\otimes V^{*}\rightarrow \C$
and $b_V':\C \rightarrow V^*\otimes V$ defined by
\begin{align*}
  d_V'&=d_Vc_{V,V^*}(\theta_V \otimes \Id_{V^*}) & b_V'&=(\Id_{V^*} \otimes
  \theta_V)C_{V,V^*}b_V.
\end{align*}  A direct computation shows that
\begin{align*}
d_V'(v\otimes f)&=f(K^{1-N}v) & b_V'(1)&=\sum K^{N-1}v_i \otimes v_i^*.
\end{align*}
\begin{lemma}\label{L:E3typ0}
 Let $V_\lambda$ be typical $\Ubar$-module then $\dim_\cat(V_\lambda)=0$.
\end{lemma}
\begin{proof}
  Let $\{v_i\}$ be a basis of $V_\lambda$ such that $v_i$ is a non-zero vector
  of weight $\lambda-2i$.  By definition we have
  $\dim_\cat(V_\lambda)=(d_{V_\lambda}'\circ b_{V_\lambda})(1)$ which is equal
  to:
 $$\sum_{i=0}^{N-1} v_i \otimes v_i^*=  \sum_{i=0}^{N-1} v_i^*(K^{1-N}v_i)=
 \sum_{i=0}^{N-1} q^{(N-1)(\lambda-2i)} =
 q^{(N-1)\lambda}\frac{1-q^{-2N}}{1-q^{-2}}$$ where $q^{-2N}=1$ and so
 $\dim_\cat(V_\lambda)=0$.
\end{proof}
\begin{lemma}\label{L:TS'}
  Let $V_\lambda$ be a typical module and let $V_{\lambda'}$ be any simple
  weight module with highest weight $\lambda'$.  We have
\begin{align*}
  S'(V_\lambda,V_{\lambda'})&=q^{(\lambda +1 -N)(\lambda' +1
    -N)}\frac{\{N(\lambda' +1 -N)\}}{\{\lambda' +1 -N\}},
\end{align*}
where $\frac{\{N(\lambda' +1 -N)\}}{\{\lambda' +1 -N\}}$ is a Laurent
polynomial in $q^{\lambda'}$.
\end{lemma}
\begin{proof}
  The proof follows from a direct computations.  A detailed presentation of an
  analogous computation is given in Proposition 2.2 of \cite{GP2}.
\end{proof}
Let $I$ be the set of isomorphism classes of simple objects of $\cat$ and let
$B$ be the subset of $ I$ consisting of isomorphism classes of typical
$\Ubar$-modules.
\begin{lemma}
  The usual invariant $F$ restricted to $\links_B$ is zero.
\end{lemma}
\begin{proof}
  Follows from Corollary~\ref{C:dim0} and Lemma \ref{L:E3typ0}.
\end{proof}
Let $V_\lambda$ be in $B$ then from Lemma \ref{L:TS'} we have
$S'(V_\lambda,V)\neq 0$ for all $V\in B$ and $S'(V_\lambda,W)= 0$ for all
$W\in I \setminus B$.  In other words, the set $B$ is complete.
\begin{prop}
  Every element of $J\in B$ is ambidextrous and $B=A(J)$.  The construction in
  Section~\ref{SS:inv} gives a function $\qd_J:B\rightarrow \C$ and an
  invariant $F'_J$, which do not depend on $J$ up to multiplication by a
  non-zero element of $\C$.
\end{prop}
 \begin{proof}
   Since $B$ is complete and contains $J$ then the lemma follows from
   Lemma~\ref{L:B}.
 \end{proof}

 Let $\qd=\qd_{J_0}$ and $F'=F'_{J_0}$ be objects defined in Section
 \ref{SS:inv} arising from the ambidextrous element $J_0$ and the constant
 $d_0=1/(\prod_{j=0}^{N-2}\qn{\lambda_0 +N -j})$. We will now compute $\qd$
 explicitly. By a direct computation, for $\lambda\in\C\setminus \Z$ one has
 $\frac{1}{\prod_{j=0}^{N-2}\qn{\lambda +N -j}} =
 (-1)^Nq^{-N(N+1)/2}\frac{\{\lambda +1 -N\} }{\{N(\lambda +1 -N)\}}$.
 Therefore, from the expression of $S'$ in Lemma \ref{L:TS'} we have that\\
 ${\prod_{j=0}^{N-2}\qn{\lambda_0 +N-j}}S'(J,V_\lambda)=
 {\prod_{j=0}^{N-2}\qn{\lambda +N -j}}S'(V_\lambda,J)$.  Then by definition we
 have
 $$
 \qd(V_\lambda)=\frac{1}{\prod_{j=0}^{N-2}\qn{\lambda +N -j}}
 $$
 where we use the above choice of $d_0$.

 Next we will show that $F'$ restricts to the colored Alexander invariant
 given by Jun Murakami in \cite{jM}.  The colored Alexander invariant is a
 reconstruction of the link invariants defined in \cite{ADO}.  Murakami's
 construction uses the universal $R$-matrix of $\Uslt$ and state sums.  In
 particular, let $T$ be a (1,1)-tangle whose $i$th component is colored by a
 parameter $\lambda_i$ ($i=1,...,k$) and the first component is the open
 component.  Murakami defines $O_T^N(\lambda_1,...,\lambda_k)$ to be the
 element of $\End(V_{\lambda_1})$ obtained by assigning the matrix elements of
 the $R$-matrix for the crossings of $T$ and particular scalars to the maximal
 and minimal points of $T$ (these scalars are same as the scalars coming from
 the morphisms $b_{V_{\lambda_i}}, d_{V_{\lambda_i}}, b'_{V_{\lambda_i}},
 d'_{V_{\lambda_i}}$ given above).  Let
 $\Phi_T^N(\lambda_1,...,\lambda_k)=\qd(V_{\lambda_1})
 O_T^N(\lambda_1,...,\lambda_k)$.  Then in \cite{jM} Murakami shows that
 $\Phi_T^N(\lambda_1,...,\lambda_k)$ is a framed version of the analogous
 invariant defined in \cite{ADO}.  Thus, $\Phi_T^N(\lambda_1,...,\lambda_k)$
 is a well defined invariant of a colored framed link $L$ obtained by closing
 the tangle $T$; denote this invariant by $\Phi_L^N$.  Since the construction
 of $F'$ uses the same $R$-matrix, duality, twist and scaling $\qd$ as the
 construction of $\Phi_L^N$ we have proved the following theorem.

\begin{theorem}
  The invariant $F'$ restricted to framed links colored with typical modules
  is equal to the colored Alexander invariant $\Phi_L^N$.
\end{theorem}

\linespread{1}

\end{document}